# ARTIFICIAL FREE CLONE OF SIMPLEX METHOD FOR FEASIBILITY

Muhammad Imtiaz[1], Nasir Touheed[2] and Syed Inayatullah[3]

## Abstract

This paper presents a method which is identical to simplex method phase 1, but do not need any artificial variable (or artificial constraints). So, the new method works in original variable space but visits the same sequence of pivots as simplex method does. Recently, (Inayatullah, Khan, Imtiaz & Khan, New artificial-free Phase 1 Simplex Method, 2009) claimed a similar method, here in this paper we have presented a counter example which shows in some special conditions of degeneracy their method may deviate from the simplex path.

Because of its simplicity, the method presented in this paper is highly classroom oriented. So, indeed there is no need to work with artificial variables (or artificial constraints) in simplex method any more.

---

[1] Department of Mathematical Sciences, University of Karachi, Karachi, Pakistan. imtiaz@uok.edu.pk
[2] Institute of Business Administration, Karachi, Pakistan. ntouheed@iba.edu.pk
[3] Department of Mathematical Sciences, University of Karachi, Karachi, Pakistan. inayat@uok.edu.pk

# 1 Introduction

The Linear Programming Problem is by far the most widely used optimization model. Its impact on economic and government modeling is remarkable. Simplex method is most useful tool to teach and solve practical linear programming problems. Since 1963, it has been the preferred method of LP practitioners (Shamir, 1987). But it often becomes inadequate and laborious to solve programs without any given initial basis. To solve a linear programming problem by simplex method the foremost need is the knowledge of a basic feasible solution. If in an LP the initial basic feasible solution is unknown, we apply the simplex method in two phases (Dantzig, Orden, & Wolfe, 1955) (Wagner, 1956) called phase 1 and phase 2. In phase 1 (non-negative) artificial variables are added to the problem to obtain a basic feasible solution artificially, resulting in an additional objective function whose value is equal to the minimization of the sum of all the artificial variables. This objective function is called phase 1 objective of the LP. In this paper, we call it as feasibility objective. The purpose of phase 1 process is to maintain the feasibility and minimize the sum of artificial variables as possible. A zero objective value at the end of phase 1 indicates that all the artificial variables have reached to value zero and hence our current basis is feasible to the original problem. We then come back to our original objective and solve it by simplex phase 2. In the other case, we conclude that the problem has no solution.

Although the approach mentioned above is the most customary but alternate approaches are also available. Zoutendijk (1976) presents different variations of the phase 1 simplex method. He has also presented an artificial free Big M like method. Arsham (1997a) and (1997b) proposed alternate but artificial-free methods to perform phase 1. Because of matured reliability and better complexity in practical problems, simplex method still preferable over all above mentioned artificial free approaches. Artificial-free methods are valuable as they do not require a basic feasible solution to start-with and they are space efficient as well. Recently Inayatullah, Khan, Imtiaz &

Khan (2009) have presented another artificial-free version of phase 1 process (say AFM), which was a modified form of the method presented in Zoutendijk (1976). They claimed that the version is geometrically analogous to the simplex phase 1 process of artificial variables, that is the order of the visited corner points is identical to the simplex phase 1. In section 2, a counterexample has been presented where the claim proved invalid. Actually their claim is true when in simplex phase 1 degeneracy do not occur due to artificial variables, but if simplex phase 1 encounters degenerate artificial variables then AFM method may deviate from the path of simplex method that is visited corner points may not be the same.

Here, in this paper we are developing another artificial-free version, which is a true clone of the simplex phase 1. It shall follow same sequence of pivoting iterations as well as the corner points. In section 6, we shall also describe the dual counterpart of our artificial free version, which is indeed an artificial constraint free version of traditional dual simplex phase 1.

## 2 Counter Example

For the following system of inequalities

$$\begin{aligned} x_1 + x_2 &\geq 4 \\ x_1 + 2x_2 &\geq 6 \\ x_1 - x_2 &\leq 2 \\ 5x_1 + 2x_2 &\geq 10 \\ x_1 \geq 0, x_2 &\geq 0 \end{aligned}$$

Path of Simplex Method phase 1: (0, 0), (2,0), (1, 5/2), (2/3, 10/3)
Path of AFM for feasibility: (0, 0), (2,0), (3, 1), (10/3, 4/3).

Here, AFM deviated from the path of simplex phase 1 that is visited corner points are not same.

## 3 Deduction of Artificial-free Form from Auxiliary Form of LP

Consider the following LPP,

(1)
$$\begin{aligned} &\textit{Maximize } Z = c^T x \\ &Ax \leq b, \\ &x \geq 0,\ x \in \Re^p \\ &A \in \Re^{m \times p},\ b \in \Re^m,\ c \in \Re^p \end{aligned}$$

Here $b$ would not necessarily be completely non-negative. By adding the unrestricted slack vector $x_B$ and denoting the original variable vector $x$ as $x_N$, we can have an equivalent equality form of the above system,

(2)
$$\begin{aligned} &Ax_N + x_B = b, \\ &x_N \geq 0,\ x_N \in \Re^p,\ x_B \in \Re^m \end{aligned}$$

Clearly, for above system the readily available basis is $B$. However, $B$ may not constitute a feasible basis (because of some negative values in $b$). We can decompose $x_B$ into difference of two non-negative variables $x_B^+$ and $x_B^-$.

(3)
$$\begin{aligned} &Ax_N + x_B^+ - x_B^- = b \\ &x_N \geq 0,\ x_B^+ \geq 0,\ x_B^- \geq 0 \\ &x_N \in \Re^p,\ x_B^+ \in \Re^m,\ x_B^- \in \Re^m \end{aligned}$$

Following is the policy for variables from $x_B^+ \cup x_B^-$ to reside or not to reside into the basis.

"For $i \in B$ when $b_i \geq 0$ then $x_i^+$ would be the basic variable and when $b_i < 0$ then $x_i^-$ would be the basic."

If $x_B^-$ is exclusive to the current basis then the current basis is feasible and we can start simplex phase 2 directly, otherwise to obtain a feasible basis we transform such system into the following linear program (often called as auxiliary form),

If $L = \{l : b_l < 0\}$, clearly $l$ is a subset of $B$.

(4)
$$\text{Minimize } \sum_L (a_{iN} x_N) + x_L^+ - \sum_L b_i$$
$$A x_N + x_B^+ - x_B^- = b$$
$$x \geq 0, x \in \Re^p, x_B \in \Re^m$$

Purpose of the objective function $\text{Minimize } \sum_L (a_{iN} x_N) + x_L^+ - \sum_L b_i$ (called the feasibility objective) is to force the sum of $x_L^-$ to zero, which ultimately (for feasible problems) makes the current basis a feasible one. By using the policy (described earlier in this section), we may construct the following pivoting table with bases $x_{B-L}^+$ and $x_L^-$.

|  | $x_N$ | $x_{B-L}^+$ | $x_{B-L}^-$ | $x_L^+$ | $x_L^-$ | $b$ |
|---|---|---|---|---|---|---|
| $w$ | $\sum_L a_{iN}$ | $\mathbf{0}_{1 \times (B-L)}$ | $\mathbf{0}_{1 \times (B-L)}$ | $\mathbf{1}_{1 \times L}$ | $\mathbf{0}_{1 \times L}$ | $\sum_L b_i$ |
| $z$ | $c_N^T$ | $\mathbf{0}_{1 \times (B-L)}$ | $\mathbf{0}_{1 \times (B-L)}$ | $\mathbf{0}_{1 \times L}$ | $\mathbf{0}_{1 \times L}$ | $0$ |
| $x_{B-L}^+$ | $A_{(B-L) \times N}$ | $I_{B-L}$ | $-I_{B-L}$ | $\mathbf{0}_{(B-L) \times L}$ | $\mathbf{0}_{(B-L) \times L}$ | $b_{B-L}$ |
| $x_L^-$ | $a_{L \times N}$ | $\mathbf{0}_{L \times (B-L)}$ | $\mathbf{0}_{L \times (B-L)}$ | $I_L$ | $-I_L$ | $b_L$ |

Here $x_B^-$ are analogous to the artificial variables in usual simplex method. Obviously, once they leave there is no reason to enter them again into the basis. So, in the pivoting table we do not need columns of coefficients of $x_B^-$. By eliminating these unnecessary columns pivoting table size would be reduced we may call it reduced pivoting table.

$$\begin{array}{c} \\ w \\ z \\ x_{B-L}^+ \\ x_L^- \end{array} \begin{array}{c} x_N \quad\quad x_{B-L}^+ \quad\quad x_L^+ \quad\quad b \\ \left[ \begin{array}{cccc|c} \sum_L a_{iN} & \mathbf{0}_{1\times(B-L)} & \mathbf{1}_{1\times L} & \sum_L b_i \\ c_N^T & \mathbf{0}_{1\times(B-L)} & \mathbf{0}_{1\times L} & 0 \\ \hdashline A_{(B-L)\times N} & I_{B-L} & \mathbf{0}_{(B-L)\times L} & b_{B-L} \\ a_{L\times N} & \mathbf{0}_{L\times(B-L)} & I_L & b_L \end{array} \right] \end{array}$$

Interestingly, if we ignore the superscripts '+' and '−' then this pivoting table would looks quite similar to a normal simplex table despite that it has some infeasible basic variables which have '+1' coefficients in phase 1 objective function $w$.

In this method, we are opting procedure to find the entering basic variable similar to simplex method i.e. we may identify the entering basic variable by seeking most negative element in $w$-row (excluding last element, which is actually equals to feasibility objective value). After determining entering basic variable the leaving variable is determined by taking following minimum-ratio test. Ratios related to variables in $x_{B-L}^+$ are obtainable by dividing $b_{B-L}$ by corresponding element in pivot column (only if that element is positive) and ratios related to variables in $x_L^-$ is obtainable by dividing $b_L$ by corresponding element in pivot column (only if that element is negative). The element in the pivot column corresponding to minimum ratio would be the pivot element. Let the pivot element is $x_{rs}$, if $r \in L$ then set $L := L - r$ and perform the pivot operation. Now in the new pivoting table proceed with the same procedure of entering and leaving basic variable until $L$ becomes empty.

## 4 Artificial-free Simplex Method in Dictionary Form

Consider the following LP problem,

(5)
$$\begin{aligned} \text{Maximize} \quad & c^T x \\ \text{subject to} \quad & Ax \leq b \\ & x \geq 0, \, x \in \Re^n \end{aligned}$$

where $x$ is the decision variable vector, $A \in \Re^{m \times n}$, $b \in \Re^m$, $c \in \Re^n$.

Given a basis $B$ we call the following matrix $D(B)$, a simplex tableau or a dictionary (Chvatal, 1983) associated with basis $B$.

$$D(B) = \begin{bmatrix} z & -c^T \\ b & A \end{bmatrix}$$

Here elements of $D(B)$ are denoted by the array $d_{IJ}$, where $I = \{0, B\}$ and $J = \{0, N\}$.

The associated primal basic solution could directly be obtained by setting $x_N = 0$ and $x_B = b$. For $i \in B$ and $j \in N$, if the component $d_{ij}$ of $D$ is non-zero, then one can obtain an equivalent dictionary with new basis $B := B - \{i\} + \{j\}$, the new non-basis $N := N + \{i\} - \{j\}$. The replacement is known as a *pivot operation* on $(i, j)$.

It is well known that if $d_{B0} \geq 0$ then $B$ (or $D(B)$) is called primal feasible; if $d_{0N} \geq 0$ then $B$ (or $D(B)$) is called dual feasible; if $B$ (or $D(B)$) is both primal and dual feasible then it is optimal. A basis $B$ (or a dictionary) is called primal inconsistent if there exists $i \in B$ such that $d_{i0} < 0$ and $d_{iN} \geq 0$, and dual inconsistent if there exists $j \in N$ such that $d_{0j} \leq 0$ and $d_{Bj} < 0$.

Consider the artificial-free form deduced in section 3, in the dictionary notation, the basic columns are omitted from the pivoting table. So, by deleting the basic columns we get the following dictionary,

$$\begin{array}{c} \\ w \\ z \\ x^+_{B-L} \\ x^-_L \end{array} \begin{array}{c} b \\ \left[\begin{array}{c} \sum_L b_i \\ 0 \\ \hline b_{B-L} \\ b_L \end{array}\right. \end{array} \begin{array}{c} x_N \\ \sum_L a_{iN} \\ c_N^T \\ \hline A_{(B-L) \times N} \\ a_{L \times N} \end{array} \begin{array}{c} x_L^+ \\ \mathbf{1}_{1 \times L} \\ \mathbf{0}_{1 \times L} \\ \hline \mathbf{0}_{(B-L) \times L} \\ I_L \end{array}\left. \right]$$

Now, here dictionary size can be more reduced by the fact that column of $x_L^+$ is identical to the hidden basic column of $x_L^-$. So, we can also eliminate column of $x_L^+$ with a condition that for every pivot on $x_{rs}$ with $r \in L$, after the pivot we must replace newly produced non-basic column of $x_r^-$ by pre-calculated column of $x_r^+$. Actually it is not difficult, technically columns of $x_r^-$ and $x_r^+$ are same with the exception that first element of $x_r^+$ is equal to the first element of $x_r^-$ plus '1'.

$$\begin{array}{c} \\ w \\ z \\ x^+_{B-L} \\ x^-_L \end{array} \begin{array}{c} b \\ \left[\begin{array}{c} \sum_L b_i \\ 0 \\ \hline b_{B-L} \\ b_L \end{array}\right. \end{array} \begin{array}{c} x_N \\ \sum_L a_{iN} \\ c_N^T \\ \hline A_{(B-L) \times N} \\ a_{L \times N} \end{array}\left. \right]$$

**Problem 1**

Given a dictionary $D(B)$, obtain primal feasibility.

**Algorithm 1: Artificial-free Simplex Method for Primal Feasibility (ASM)**

**Step 1:** Let $L$ be a maximal subset of $B$ such that $L = \{i : d_{i0} < 0, i \in B\}$. If $L = \varphi$ then $D(B)$ is primal feasible. **Exit.**

**Step 2:** Denote the basic variables $x_L$ by $x_L^-$ and compute feasibility objective vector $w(B) \in \Re^N$ such that $w(B)_j = \sum_{l \in L} d_{ij}$, $j \in N$. Append $w$ to the top of the dictionary $D(B)$.

**Step 3:** Let $K \subseteq N$ such that $K = \{j : w(B)_j < 0, j \in N\}$. If $K = \varphi$ then $D(B)$ is primal infeasible. **Exit.**

**Step 4:** Choose $m \in K$ such that $w(B)_m \leq w(B)_k$, $\forall k \in K$
(Ties should be broken arbitrarily)

**Step 5:** Choose $r_1 \in B - L$ such that
$r_1 = \arg\min\{\{d_{i0} / d_{im} \mid (d_{i0} \geq 0, d_{im} > 0)\}, i \in B - L\}$
Choose $r_2 \in L$ such that $r_2 = \arg\min\{\{d_{i0} / d_{im} \mid (d_{i0} \leq 0, d_{im} < 0)\}, i \in L\}$
Set $r := \text{augmin}\left\{\dfrac{d_{r_1 0}}{d_{r_1 m}}, \dfrac{d_{r_2 0}}{d_{r_2 m}}\right\}$

**Step 6:** Make a pivot on $(r, m)$ ($\Rightarrow$ Set $B := B + \{m\} - \{r\}$, $N := N - \{m\} + \{r\}$ and update $D(B)$).

**Step 7:** If $r \in L$, $d_{0m} := d_{0m} + 1$ & notation of $x_r^-$ would replaced by $x_r^+$. Set $L := L - \{r\}$.

**Step 8:** Go to Step 3.

## 5  Proof of Correctness

Our artificial-free method (ASM) could start with an infeasible basis without making it artificially feasible. As stated earlier the objective of the simplex phase 1 is, to minimize the sum of all artificial variables. In an analogous sense, as shown in section 4 the feasibility objective $w$ is to minimize $\sum_L x_i^-$ which is equivalent to sum of rows corresponding to infeasible basic variables. Just like simplex method, our method intends to achieve the feasibility of the infeasible (negative) variables provided the feasibility of existing feasible (non-negative) variables remains preserved.

The overall structure of pivoting forces $\sum_{i \in L} b_i$ to be strictly decreasing throughout the iterations for-degenerate pivots and constant for degenerate pivots. So, finiteness of total number of bases in every LP problem proves finiteness of our method for complete non-degenerate LP problems.

**Example 1**

Obtain a feasible basis of the following LP using ASM

$$x_1 + x_2 \geq 4$$
$$x_1 + 2x_2 \geq 6$$
$$x_1 - x_2 \leq 2$$
$$5x_1 + 2x_2 \geq 10$$
$$x_1 \geq 0, x_2 \geq 0$$

By adding unrestricted slack variables $x_3$, $x_4$, $x_5$ and $x_6$, we can construct the associated dictionary along with feasibility objective vector $w$ of the above problem as

Initial table:

|       | $b$  | $x_1$ | $x_2$ |
|-------|------|-------|-------|
| $w$   | -20  | -7    | -5    |
| $x_3$ | -4   | -1    | -1    |
| $x_4$ | -6   | -1    | -2    |
| $x_5$ | 2    | 1     | -1    |
| $x_6$ | -10  | -5    | -2    |

Here $L = \{3,4,7\}$, replace $x_L \to x_L^-$

|         | $b$  | $x_1$ | $x_2$ |
|---------|------|-------|-------|
| $w$     | -20  | -7    | -5    |
| $x_3^-$ | -4   | -1    | -1    |
| $x_4^-$ | -6   | -1    | -2    |
| $x_5$   | 2    | 1*    | -1    |
| $x_6^-$ | -10  | -5    | -2    |

Iteration 1:

Here $r = 5 \in B - L$,

|       | $b$ | $x_5$ | $x_2$ |
|-------|-----|-------|-------|
| $w$   | -6  | 7     | -12   |
| $x_3^-$ | -2 | 1   | -2    |
| $x_4^-$ | -4 | 1   | -3    |
| $x_1$ | 2   | 1     | -1    |
| $x_6^-$ | 0  | 5    | -7*   |

Iteration 2:

Here $r = 6 \in L$, replace $x_6^- \to x_6^+$, $L = \{3,4,6\} - \{6\} = \{3,4\}$, $d_{02} := d_{02} + 1$

|       | $b$ | $x_5$ | $x_6^+$ |
|-------|-----|-------|---------|
| $w$   | -6  | -11/7 | -5/7    |
| $x_3^-$ | -2 | -3/7 | -2/7    |
| $x_4^-$ | -4 | -8/7* | -3/7   |
| $x_1$ | 2   | 2/7   | -1/7    |
| $x_2$ | 0   | -5/7  | -1/7    |

Iteration 3:

Here $r = 4 \in L$, replace $x_4^- \to x_4^+$, $L = \{3,4\} - \{4\} = \{3\}$, $d_{01} := d_{01} + 1$

|       | $b$ | $x_4^+$ | $x_6^+$ |
|-------|-----|---------|---------|
| $w$   | -1/2 | -3/8   | -1/8    |
| $x_3^-$ | -1/2 | -3/8* | -1/8   |
| $x_5$ | 7/2  | -7/8   | 3/8     |
| $x_1$ | 1    | ¼      | -1/4    |
| $x_2$ | 5/2  | -5/8   | 1/8     |

Iteration 4:

Here $r = 3 \in L$, replace $x_3^- \to x_3^+$, $L = \{3\} - \{3\} = \phi$,

|  | b | $x_3^+$ | $x_6^+$ |
|---|---|---|---|
| $x_4^+$ | 4/3 | -8/3 | 1/3 |
| $x_5$ | 14/3 | -7/3 | 2/3 |
| $x_1$ | 2/3 | 2/3 | -1/3 |
| $x_2$ | 10/3 | -5/3 | 1/3 |

Primal feasibility is achieved, the feasible solution is (2/3, 10/3) and the sequence of corner points are (0, 0), (2, 0), (1, 5/2), (2/3, 10/3) i.e. the corner points are same as in simplex method.

# 6   Dual Version of Artificial-free Simplex Method

**Problem 2**
Given a dictionary $D(B)$, obtain dual feasibility.

**Algorithm 2: Artificial-free Simplex Method for Dual Feasibility (ASMD)**

**Step 1:**   Let K be a maximal subset of $B$ such that $K = \{j : d_{0j} < 0, j \in N\}$. If $K = \varphi$ then $D(B)$ is primal feasible. **Exit.**

**Step 2:**   Denote the dual slack variables $y_K$ by $y_{\bar{K}}$ and compute feasibility objective vector $w'(B) \in \Re^B$ such that $w'(B)_i = -\sum_{k \in K} d_{ik}$, $i \in B$. Append $w'$ to the right of the dictionary $D(B)$.

**Step 3:** Let $L \subseteq N$ such that $L = \{i : w'(B)_i < 0, i \in B\}$. If $L = \varphi$ then $D(B)$ is primal infeasible. **Exit.**

**Step 4:** Choose $r \in L$ such that $w'(B)_r \leq w'(B)_l$, $\forall l \in L$

(Ties should be broken arbitrarily)

**Step 5:** Choose $m_1 \in B - K$ such that

$$m_1 = \arg\max\{\{d_{0j}/d_{rj} \mid (d_{0j} \leq 0, d_{rj} > 0)\}, j \in B - K\}$$

Choose $m_2 \in K$ such that $m_2 = \arg\max\{\{d_{0j}/d_{rj} \mid (d_{0j} \geq 0, d_{rj} < 0)\}, j \in K\}$

Set $m := \text{augmax}\left\{\dfrac{d_{0m_1}}{d_{rm_1}}, \dfrac{d_{0m_2}}{d_{rm_2}}\right\}$

**Step 6:** Make a pivot on $(r, m)$ ($\Rightarrow$ Set $B := B + \{m\} - \{r\}$, $N := N - \{m\} + \{r\}$ and update $D(B)$).

**Step 7:** If $m \in K$, $d_{r(|N|+1)} := d_{r(|N|+1)} + 1$ & notation of $y_m^-$ would replaced by $y_m^+$. Set $K := K - \{m\}$.

**Step 8:** Go to Step 3.

**Example 2**

Obtain a complementary dual feasible basis of the following LP by using ASMD.

$$\text{Maximize} \quad z = 5x_1 - 2x_2 + 7x_3$$
$$\text{subject to} \quad -8x_1 + x_2 - 5x_3 \geq 6$$
$$8x_1 + 8x_2 \geq 8$$
$$-7x_1 + 9x_2 \geq 1$$
$$-7x_1 + 7x_2 - 9x_3 \geq 6$$
$$x_1 \geq 0, x_2 \geq 0, x_3 \geq 0$$

By adding unrestricted slack variables $x_4$, $x_5$, $x_6$ and $x_7$ we can construct the associated dictionary along with column vector $w'$ (negative sum of columns of infeasible dual basic variables) of the above problem. Here we demonstrating explicitly the dual variables $y_1$, $y_2$, $y_3$, $y_4$, $y_5$ and $y_6$ as we required to observe dual slack variables as well

Initial table:

|       | b  | $x_1$ | $x_2$ | $x_3$ | $w'$ |       |
|-------|----|-------|-------|-------|------|-------|
| z     | 0  | −5    | 2     | −7    | 12   |       |
| $x_4$ | −6 | 8     | −1    | 5     | −13  | $y_4$ |
| $x_5$ | −8 | −8    | −8    | 0     | 8    | $y_5$ |
| $x_6$ | −1 | 7     | −9    | 0     | −7   | $y_6$ |
| $x_7$ | −6 | 7     | −7*   | 9     | −16  | $y_7$ |
|       |    | $\overline{y_1}$ | $y_2$ | $\overline{y_3}$ | |

Iteration 1:

|       | b     | $x_1$ | $x_7$ | $x_3$  | $w'$  |       |
|-------|-------|-------|-------|--------|-------|-------|
| z     | -12/7 | −3    | 2/7   | -31/7  | 52/7  |       |
| $x_4$ | -36/7 | 7*    | -1/7  | 26/7   | -75/7 | $y_4$ |
| $x_5$ | -8/7  | −16   | -8/7  | -72/7  | 184/7 | $y_5$ |
| $x_6$ | -47/7 | −2    | -9/7  | -81/7  | 95/7  | $y_6$ |
| $x_2$ | 6/7   | −1    | -1/7  | -9/7   | 16/7  | $y_2$ |
|       |       | $\overline{y_1}$ | $y_7$ | $\overline{y_3}$ | |

Iteration 2:

|       | b       | $x_4$  | $x_7$  | $x_3$    | $w'$    |         |
|-------|---------|--------|--------|----------|---------|---------|
| z     | -192/49 | 3/7    | 11/49  | -139/49  | 139/49  |         |
| $x_1$ | -36/49  | 1/7    | -1/49  | 26/49    | -26/49  | $y_1^+$ |

| | | | | | |
|---|---|---|---|---|---|
| $x_5$ | -632/49 | 16/7 | -72/49 | -88/49 | 88/49 | $y_5$ |
| $x_6$ | 257/49 | 2/7 | -65/49 | -515/49 | 515/49 | $y_6$ |
| $x_2$ | 6/49 | 1/7 | -8/49 | -37/49 | 37/49 | $y_2$ |
| | $y_4$ | $y_7$ | $y_3^-$ | | |

Iteration 3:

| | b | $x_4$ | $x_7$ | $x_1$ | |
|---|---|---|---|---|---|
| z | -102/13 | 31/26 | 3/26 | 139/26 | |
| $x_3$ | -18/13 | 7/26 | -1/26 | 49/26 | $y_3^+$ |
| $x_5$ | -200/13 | 36/13 | -20/13 | 44/13 | $y_5$ |
| $x_6$ | -121/13 | 81/26 | -45/26 | 515/26 | $y_6$ |
| $x_2$ | -12/13 | 9/26 | -5/26 | 37/26 | $y_2$ |
| | $y_4$ | $y_7$ | $y_1^+$ | | |

Dual feasibility is achieved and (0, −12/13, −18/13) be the complementary feasible solution.

## 7 Conclusion

In this paper, we developed an artificial variable free auxiliary form of simplex phase 1. On the basis of that a new method is presented "Artificial-free simplex method (ASM)" which is indeed an artificial variable free clone of simplex method. Advantages of the new method are evident from its name, that is, it do not need artificial variables so it could start with any feasible or infeasible basis of any LP. Also, it is space efficient because its dictionary do not need additional columns of artificial variables. And last but not least, its pivoting sequence is identical to simplex method so it has the same complexity as the simplex method has. Therefore, in the classroom teaching there is no need to teach the odd-looking concept of artificial variables and artificial constraints.